\documentclass[11pt,fleqn,leqno]{article}  
\usepackage{amsmath} 
\usepackage{amssymb} 
\newtheorem{thm}{Theorem}[section]
\newcommand{\ba}{\begin{array}}
\newcommand{\ea}{\end{array}} 
\newcommand{\pri}{\prime}

\newcommand{\+}{\!+\!}
\newcommand{\m}{\!-\!}

\newcommand{\z}{\zeta}
\newcommand{\cq}{\frac{v'(z)-v'(\z)}{z-\z}}
\newcommand{\SA}{{\mathbb A}}
\newcommand{\SB}{{\mathbb B}}
\newcommand{\SC}{{\mathbb C}}
\newcommand{\SD}{{\mathbb D}}
\renewcommand{\a}{\alpha}
\renewcommand{\b}{\beta}

\newcommand{\f}{\phi}

\renewcommand{\k}{\kappa}

\newcommand{\bea}{\begin{eqnarray}}
\newcommand{\eea}{\end{eqnarray}}

\newcommand{\quar}{{\lower0.25ex\hbox{\raise.6ex\hbox{\the\scriptfont0 1}
\kern-.2em\slash\kern-.1em\lower.25ex\hbox{\the\scriptfont0 4}}}}

\begin{document}
\newtheorem{lem}[thm]{Lemma}
\newtheorem{cor}[thm]{Corollary}

\title{\bf Discriminants and Functional 
Equations for  Polynomials Orthogonal on the Unit Circle }
\author{
Mourad E.H. Ismail${}^1$\thanks{Research partially supported by NSF grant DMS-99-70865}\;
and Nicholas S. Witte${}^2$}

\maketitle

\begin{abstract}
We derive raising and lowering operators for orthogonal polynomials on 
the unit circle and find second order differential and $q$-difference 
equations for these polynomials. 
A general functional equation is found which allows one to relate the
zeros of the orthogonal polynomials to the stationary values of an explicit 
quasi-energy and implies recurrences on the orthogonal polynomial 
coefficients.
We also evaluate the discriminants and quantized discriminants of polynomials 
orthogonal on the unit circle. 
\end{abstract}

\bigskip

{\bf Running Title}: Discriminants and Functional  Equations



\bigskip


\setcounter{section}{1}
\setcounter{equation}{0}
\setcounter{thm}{0}

\noindent{\bf 1. Introduction}.  
Let $w(z)$ be  a weight function  supported on a subset of the unit circle
and assume that $w$ is normalized by
\bea
 \int_{|\z|=1}  w(\z) \frac{d\z}{i\z} = 1.
\label{uc-norm}
\eea
Let $\f_n(z)$ be the  polynomials  orthonormal with respect to $w(z)$, 
that is 
\bea
 \int_{|\z|=1} \f_m(\z)\overline{\f_n(\z)} w(\z) \frac{d\z}{i\z} 
  = \delta_{m,n}.
\label{uc-ortho}
\eea
A general background to orthogonal polynomial systems defined on the unit
circle can be found in the monographs \cite{Sz2}, \cite{Ge1}, and
\cite{Fr}, while more recent surveys are to be found in \cite{Ge2}, \cite{Ge3}
and \cite{Ne} and from an interesting perspective, in the course notes of 
\cite{Ma}.
In this work we first derive raising and lowering operators for $\f_n$ 
under certain smoothness conditions on the weight function then use these 
first order operators  
to derive a linear second order differential equation satisfied by the 
orthogonal polynomials. This will be done in \S 2.  These results are unit 
circle analogues of the results of Bauldry \cite{Ba}, Bonan and Clark
\cite{Bo:Cl}, and Chen and Ismail \cite{Ch:Is}. The  
external field \cite{Is2}, \cite{Sa:To} is the function $v$ defined by 
\bea
w(z) =  \exp(-v(z)).
\label{uc-logw}
\eea
We illustrate these general results by three examples - the circular Jacobi
polynomials, the Szeg\"o polynomials and the orthogonal polynomial system
defined by the modified Bessel function.

Flowing from the results in \S 2 we derive a functional equation and relate
this to the zeros of orthogonal polynomials defined on the unit circle in
\S 3. This is the analogue of the electrostatic interpretation of the zeros
of orthogonal polynomials defined on the real line, but in this case analyticity
means that the quasi-energy function derived has stationary points at the zeros which
are saddle-points, not minima. This functional equation implies a general
relationship on the orthogonal polynomial system, which is usually expressed as a
recurrence relation on the polynomial coefficients.

In \S 4 we derive $q$-analogues of \S 2. 
The  external field is now the function $u$ defined through 
\bea
(D_{q} w)(z) = -u(qz) w(qz), 
\label{qop-logw}
\eea
where $D_q$ is the $q$-difference operator 
\bea
(D_{q} f)(z) := \frac{f(z)-f(qz)}{z-qz}.
\label{q-diff}
\eea
In other words 
\bea
w(z) = w(qz)[1-(1\m q)z u(qz)], \quad |z| = 1.
\label{ucq-w}
\eea

Recall that the discriminant $D(f_n)$ of a polynomials $f_n$ is defined 
by \cite{Di}
\bea
D(f_n) = \gamma^{2n-2}\prod_{1 \le j < k \le n} (z_j-z_k)^2, 
\quad {\rm if} \quad f_n(z) = \gamma \prod_{j=1}^n (z-z_j).
\label{poly-disc}
\eea
Stieltjes \cite{St1}, \cite{St2} and Hilbert \cite{Hi}  evaluated   
the discriminants of the classical orthogonal polynomials of Hermite, 
Laguerre, and Jacobi. Schur \cite{Sc}  gave an 
interesting lemma about general orthogonal polynomials which implies 
the Stieltjes-Hilbert results when applied 
to  the  Hermite, Laguerre, and Jacobi polynomials. In \S 4 we prove 
an analogue of Schur's lemma for polynomials orthogonal on the unit circle 
and use it to give a general theorem on the evaluation of discriminants of 
orthogonal polynomials on the unit circle. This is the unit circle analogue 
of the results in \cite{Is1}. It was observed in 
\cite{Is4} that, in general,  the discriminant (\ref{poly-disc}) of
 $q$-orthogonal polynomials does not have a closed form. The appropriate 
 discriminant for discrete $q$-orthogonal polynomials is 
\bea
D(f_n, q) = \gamma^{2n-2}\, q^{n(n-1)/2}\, \prod_{1 \le j < k \le n} 
(q^{1/2}z_j - z_kq^{-1/2}) (q^{-1/2}z_j - z_kq^{1/2}), 
\label{qop-disc}
\eea
if $f_n$ is as in (\ref{poly-disc}). The above discriminant also has the 
alternate representation 
\bea
D(f_n, q) = \gamma^{2n-2}\, q^{n(n-1)/2}\, \prod_{1 \le j < k \le n} 
\left[ z_j^2+z_k^2 - z_jz_k (q+q^{-1})\right].  
\label{qop-disc2}
\eea
In particular for a quadratic polynomial $Az^2 + Bz +C$ the $q$-discriminant 
is $qB^2 - (1+q)^2 AC.$  

In \S 5 we give an expression for the $q$-discriminant of polynomials 
orthogonal on the unit circle in terms of the coefficients in the recurrence 
relations satisfied by the polynomials. As an illustration we evaluate 
the $q$-discriminant of the Rogers-Szeg\H{o} polynomials \cite{Sz1}.   

\bigskip

\setcounter{section}{2}
\setcounter{thm}{0}
\setcounter{equation}{0}

\noindent{\bf 2. Differential Equations}.
Recall that if $f$ is a polynomial of degree $n$ then the reciprocal polynomial is
\bea
f^*(z) := \sum_{k=0}^n \overline{a_k}\, z^{n-k}, \quad {\rm if} \quad 
f(z) = \sum_{k=0}^n a_k \, z^{k}, \quad {\rm and}\; a_n \ne 0,   
\label{uc-conj}
\eea
\cite{Sz2}. Let  $\f_n$ satisfy (\ref{uc-ortho}) and 
\bea
\f_n(z) = \k_n z^n + l_nz^{n-1} + {\rm lower \; order \; terms}, \quad 
 \k_n > 0 \; {\rm for} \; n > 0.
\label{uc-coeff}
\eea
Then the $\f_n$'s satisfy the recurrence relations \cite[(11.4.6), (11.4.7)]
{Sz2} 
\bea
\k_n z \f_n(z) & = & \k_{n+1} \f_{n+1}(z) - \f_{n+1}(0) \f_{n+1}^*(z),
\label{uc-rec1} \\
\k_n  \f_{n+1}(z) & = & \k_{n+1} z \f_{n}(z) + \f_{n+1}(0) \f_{n}^*(z).
\label{uc-rec2}
\eea
If we eliminate $\f_n^*$ between (\ref{uc-rec1}) and (\ref{uc-rec2}) we get the 
three term recurrence relation (XI.4, p. 91 in \cite{Ge3})
 \bea
 \k_n\f_n(0)\f_{n+1}(z) + \k_{n-1}\f_{n+1}(0)z\f_{n-1}(z) 
 =\left[ \k_n \f_{n+1}(0) + \k_{n+1}\f_n(0)z \right]\, \f_n(z). 
\label{uc-3term}
 \eea

The $\k$'s and $\f_n(0)$ are related through \cite[(11.3.6)]{Sz2}
\bea
\k_n^2 = \sum_{k=0}^n |\f_k(0)|^2.
\label{uc-k}
\eea
Thus $\k_n$ ($>0$) can be found from the knowledge of $|\f_k(0)|$. 
By equating coefficients of $z^n$ in (\ref{uc-3term}) and in view of 
(\ref{uc-coeff}) we find
\bea
  \k_n l_{n+1}\f_n(0) + \k^2_{n-1}\f_{n+1}(0) 
  = \k_n^2\f_{n+1}(0) + \k_{n+1}l_n\f_n(0). \nonumber
\label{uc-kl1}
\eea
Thus 
\bea
\k_n l_{n+1} = \k_{n+1} l_n + \overline{\f_n(0)}\f_{n+1}(0).
\label{uc-kl}
\eea
From (\ref{uc-kl}) it is possible to express $l_n$ in terms of the $\k$'s and 
$\f_j(0)$'s
\bea
  l_n = \k_n \sum^{n-1}_{j=0}
        { \overline{\f_j(0)}\f_{j+1}(0) \over \k_j\k_{j+1} }.
\label{uc-l}
\eea
The analogue of the Christoffel-Darboux formula is 
\bea
\sum_{k=0}^n \overline{\f_k(a)} \f_k(z) = 
\frac{\overline{\f_{n+1}^*(a)} \f_{n+1}^*(z) - 
\overline{\f_{n+1}(a)} \f_{n+1}(z)}{1- \overline{a}\, z}. 
\label{uc-cd}
\eea
 
\begin{thm} 
Let $w(z)$ be differentiable in a neighborhood of the unit circle, 
has moments of all integral orders and assume that the integrals 
\bea
\int_{|\z| =1} \cq \z^n w(\z) \frac{d\z}{i\z} \nonumber 
\label{uc-int}
\eea
exist for all integers $n$. Then the corresponding orthonormal polynomials 
satisfy the differential relation 
\bea
\f_n^{\pri}(z) &=& n \frac{\k_{n-1}}{\k_n} \f_{n-1}(z) 
-i \f_n^*(z) \int_{|\z| =1} \cq \f_n(\z) \overline{\f_n^*(\z)}w(\z) d\z 
 \label{uc-diffr} \\
&{}& \quad +i \f_n(z) \int_{|\z| =1} \cq \f_n(\z) \overline{\f_n(\z)}w(\z) 
d\z \nonumber 
\eea
\end{thm}
{\bf Proof.} Using the orthogonality relation (\ref{uc-ortho}) we express 
$\f_n'(z)$ in terms of the $\f_k$'s as 
\bea
\f_n^{\pri}(z) &=&  \sum_{k=0}^{n-1} \f_k(z) \, \int_{|\z| = 1}  
\f_n^{\pri}(\z) 
\overline{\f_k(\z)} w(\z) \frac{d\z}{i\z}  \nonumber\\ 
&=& \sum_{k=0}^{n-1} \f_k(z) \, \int_{|\z| = 1} 
\left[v'(\z) \overline{\f_k(\z)} + \overline{\z \f_k(\z)} 
+ \overline{\z^2\f_k^{\pri}(\z)}\right] \f_n(\z)w(\z) \frac{d\z}{i\z}, 
  \nonumber 
\label{uc-diff1}
\eea
where we have integrated by parts, then rewritten the derivative of the 
conjugated polynomial in the following way
\bea
   \frac{d}{d\z}\overline{\f_n(\z)} = -\overline{\z}^2\overline{\f'_n(\z)},
\label{uc-pf1}
\eea
and used the fact $\overline{\z} = 1/\z$ for $|\z| =1$. 
Now the orthogonality relation (\ref{uc-ortho}) and (\ref{uc-coeff}) give 
\bea
\f_n^{\pri}(z) &=& \int_{|\z| = 1} v'(\z) \f_n(\z) 
\sum_{k=0}^{n-1} \overline{\f_k(\z)}\f_k(z) \, w(\z)\frac{d\z}{i\z}\nonumber\\ 
&{}& \; + \f_{n-1}(z)  \int_{|\z| = 1}\left[ \overline{\z \f_{n-1}(\z)} 
+ \overline{\z^2\f_{n-1}^{\pri}(\z)}\right] \f_n(\z)w(\z) \frac{d\z}{i\z} 
 \nonumber\\ 
  &=& \int_{|\z| = 1}  \cq \f_n(\z) \left[\overline{\f_n^*(\z)}
\f_n^*(z) -\overline{\f_n(\z)} \f_n(z)\right] w(\z)\frac{d\z}{i}  \nonumber\\
&{}& \; + \f_{n-1}(z) \left[\frac{\k_{n-1}}{\k_n} + (n-1) 
\frac{\k_{n-1}}{\k_n}\right].   \nonumber
\label{uc-diff2}
\eea
This establishes (\ref{uc-diffr}) and completes the proof. 

We next apply (\ref{uc-rec1}) to eliminate $\f_n^*$ from (\ref{uc-diffr}),
assuming $ \f_n(0) \neq 0 $.
The result is 
\bea
\f_n^{\pri}(z) &=& n \frac{\k_{n-1}}{\k_n} \f_{n-1}(z) 
+i \frac{\k_{n-1}}{\f_n(0)} z \f_{n-1}(z) 
\int_{|\z| =1} \cq \f_n(\z) \overline{\f_n^*(\z)}w(\z) d\z \label{uc-diff3} \\
&{}& \quad +i \f_n(z) \int_{|\z| =1} \cq \f_n(\z) 
\left[\overline{\f_n(\z)} - \frac{\k_n}{\f_n(0)}\overline{\f_n^*(\z)}
\right] w(\z) d\z. \nonumber
\eea
Observe that $\overline{\f_n(\z)} - \frac{\k_n}{\f_n(0)}\overline{\f_n^*(\z)}$ 
is a polynomial of degree $n-1$. 

Let 
\bea
A_n(z) &=& n \frac{\k_{n-1}}{\k_n} + i \frac{\k_{n-1}}{\f_n(0)} z  
\int_{|\z| =1} \cq \f_n(\z) \overline{\f_n^*(\z)}w(\z) d\z,
\label{uc-A}  \\
B_n(z) &=& -i \int_{|\z| =1} \cq \f_n(\z) 
\left[\overline{\f_n(\z)} - \frac{\k_n}{\f_n(0)}\overline{\f_n^*(\z)}
\right] w(\z) d\z.
\label{uc-B}
\eea
For future reference we note that $ A_0 = B_0 = 0 $ and
\bea
   A_1(z) & = &
   \k_1-\f_1(z)v'(z)-{\f_1^2(z) \over \f_1(0)}M_1(z),
   \\
   B_1(z) & = & 
   -v'(z)-{\f_1(z) \over \f_1(0)}M_1(z),
\label{uc-lowB}
\eea
where the first moment $ M_1 $ is defined by
\bea
   M_1(z) = \int_{|\z|=1} \z \cq w(\z){d\z \over i\z}.
\label{uc-M1}
\eea

Now rewrite (\ref{uc-diff3}) in the form 
\bea
\f_n^{\pri}(z) = A_n(z)\f_{n-1}(z) - B_n(z)\f_{n}(z).
\label{uc-diff}
\eea

Define differential operators $L_{n,1}$ and $L_{n,2}$ by 
\bea
L_{n,1} &=& \frac{d}{dz} + B_n(z),  
\label{uc-L1} \\
L_{n,2} &=& -\frac{d}{dz} - B_{n-1}(z) + \frac{A_{n-1}(z)\k_{n-1}}{z\k_{n-2}} 
 + \frac{A_{n-1}(z)\k_{n }\f_{n-1}(0)}{\k_{n-2}\f_n(0)}.
\label{uc-L2}
\eea
After the elimination of $\f_{n-1}$ between (\ref{uc-diff3}) and 
(\ref{uc-3term}) we find that 
the operators $L_{n,1}$ and $L_{n,2}$ are annihilation and creation operators 
in the sense that they satisfy 
\bea
L_{n,1} \f_n(z) = A_n(z) \f_{n-1}(z), \quad L_{n,2} \f_{n-1}(z) = 
\frac{A_{n-1}(z)}{z}\, \frac{\f_{n-1}(0)\k_{n-1}}{\f_n(0)\k_{n-2}}\f_n(z). 
\label{uc-ac}
\eea
Hence we have established the second order differential equation
\bea
L_{n,2}\left(\frac{1}{A_n(z)}L_{n,1}\right) \f_n(z) =
 \frac{A_{n-1}(z)}{z}\, \frac{\f_{n-1}(0)\k_{n-1}}{\f_n(0)\k_{n-2}}\f_n(z), 
\label{uc-2ode1}
\eea 
which will also be written in the following way
\bea
  \f_n'' + P(z)\f_n' + Q(z)\f_n = 0.
\label{uc-2odePQ}
\eea

It is worth mentioning that, unlike for polynomials orthogonal on the line, 
$L_{n,1}^*$ is not related to $L_{n,2}$. In fact if we let 
\bea
(f,g) := \int_{|\z| =1} f(\z)\, \overline{g(\z)}\, w(\z)   \frac{d\z}{\z},
\label{uc-ip}
\eea 
then in the Hilbert space endowed with this inner product, the adjoint of 
$L_{n,1}$ is
\bea
(L_{n,1}^* f)(z) = z^2 f^\pri(z) + z f(z)+ [\overline{v(z) + B_n(z)}] f(z).
\label{uc-adL1}
\eea
To see this use integration by parts and the fact that for $|\z| =1$, 
$\overline{g(\z)} = \overline{g}(1/\z)$.

\noindent{\bf Example 1}
The circular Jacobi orthogonal polynomials (CJ) are defined with respect
to the weight function 
\bea
  w(z) = {\Gamma^2(a\+ 1) \over 2\pi\Gamma(2a\+ 1)}|1-z|^{2a}
\label{cj-weight}
\eea
for real $ a $ appropriately restricted. 
We find these to be classical in the sense of being related to
classical orthogonal polynomials defined on the real line and therefore
possessing their properties. They arise in a class of random unitary matrix
ensembles, the CUE, where the parameter $ a $ is related to the charge of
an impurity fixed at $ z = 1 $ in a system of unit charges located on the
unit circle at the complex values given by the eigenvalues of a member of
this matrix ensemble \cite{Wi:Forr}.
The orthonormal polynomials are
\bea
    \f_n(z) = {(a)_n \over \sqrt{n!(2a\+ 1)_n}}
    \,{}_2F_1(-n,a\+ 1;-n\+ 1\m a;z)\ ,
\label{cj-poly}
\eea
and the coefficients are
\bea
  \k_n    & = & {(a\+ 1)_n \over \sqrt{n!(2a\+ 1)_n}}
  \quad n \geq 0,\\
  l_n     & = & {na \over n\+ a}\k_n
  \quad n \geq 1,\\
  \f_n(0) & = & {a \over n\+ a}\k_n
  \quad n \geq 0.
\label{cj-klphi}
\eea
The reciprocal polynomials are
\bea
    \f^*_n(z) = {(a\+ 1)_n \over \sqrt{n!(2a\+ 1)_n}}
    \,{}_2F_1(-n,a;-n\m a;z).
\label{cj-conj}
\eea
Using the differentiation formula and some contiguous relations for 
the hypergeometric functions, combined in the form
\begin{multline}
   (1\m z){d \over dz}{}_2F_1(-n,a\+ 1;1\m n\m a;z) =
   {n(n\+ 2a) \over n\m 1\+ a} \,{}_2F_1(1\m n,a\+ 1;2\m n\m a;z)
   \\ - n \,{}_2F_1(-n,a\+ 1;1\m n\m a;z),
\label{hyp-diff}
\end{multline}
one finds the differential-recurrence relation
\bea
   (1-z)\f'_n = -n\f_n + [n(n\+ 2a)]^{1/2}\f_{n-1},
\label{cj-diff}
\eea
and the coefficient functions
\bea
  A_n(z) & = & { \sqrt{n(n\+ 2a)} \over 1-z },\\
  B_n(z) & = & { n \over 1-z }.
\label{cj-AB}
\eea
The second order differential equation becomes
\bea
  \f''_n + \f'_n \left\{ {1\m n\m a \over z} - {2a\+ 1 \over 1-z} \right\}
  + \f_n {n(a\+ 1) \over z(1\m z)} = 0.
\label{cj-2ode}
\eea

\noindent{\bf Example 2}. 
We consider a generalization of the previous example, to the situation where
\bea
  w(z) = 2^{-1-2a-2b}{\Gamma(a\+ b\+ 1) \over \Gamma(a\+ 1/2)\Gamma(b\+ 1/2)}
         |1-z|^{2a}|1+z|^{2b},
\label{sz-w}
\eea
with $ x=\cos\theta $, and the associated orthogonal polynomials are known
as Szeg\"o polynomials \cite{Sz2}. They are related to the Jacobi polynomials 
via the projective mapping of the unit circle onto the interval $ [-1,1] $,
$ z \mapsto  \dfrac{1}{2}[z\+z^{-1}] = x = \cos\theta $,
\bea
 z^{-n} \f_{2n}(z) = \SA P^{(a-1/2,b-1/2)}_{n}(\dfrac{1}{2}[z\+z^{-1}])
       +\dfrac{1}{2} \SB [z\m z^{-1}]P^{(a+1/2,b+1/2)}_{n-1}(\dfrac{1}{2}[z\+z^{-1}]),
\label{sz-poly1} \\
 z^{1-n} \f_{2n-1}(z) = \SC P^{(a-1/2,b-1/2)}_{n}(\dfrac{1}{2}[z\+z^{-1}])
       +\dfrac{1}{2} \SD [z\m z^{-1}]P^{(a+1/2,b+1/2)}_{n-1}(\dfrac{1}{2}[z\+z^{-1}]).
\label{sz-poly2}
\eea
In their study of the equilibrium positions of charges confined to the unit circle subject
to logarithmic repulsion
Forrester and Rogers considered orthogonal polynomials defined on $ x $ which are just
the first term of \ref{sz-poly1}.
Using the normalization amongst the even and odd sequences of polynomials,
orthogonality between these two sequences and the requirement that the
coefficient of $ z^{-n} $ on the right-hand side of (\ref{sz-poly2}) must
vanish, one finds explicitly that the coefficients are
\bea
 \SA & = & \left\{ {n!(a\+ b\+ 1)_n \over
                    (a\+ 1/2)_n(b\+ 1/2)_n} \right\}^{1/2},
\label{sz-coeffA} \\
 \SB & = & \dfrac{1}{2} \SA,
\label{sz-coeffB} \\
 \SC & = & n\left\{ {(n\m 1)!(a\+ b\+ 1)_{n-1} \over
                     (a\+ 1/2)_n(b\+ 1/2)_n} \right\}^{1/2},
\label{sz-coeffC} \\
 \SD & = & {n\+ a\+ b \over 2n} \SC.
\label{sz-coeffD}
\eea
Furthermore the following coefficients of the polynomials are found to be
\bea
 \k_{2n}   & = &  2^{-2n} { (a\+ b\+ 1)_{2n} \over
            \sqrt{n!(a\+ b\+ 1)_n(a\+ 1/2)_n(b\+ 1/2)_n} },
\label{sz-ke} \\
 \k_{2n-1} & = & 2^{1-2n} {(a\+ b\+ 1)_{2n-1} \over
              \sqrt{(n\m 1)!(a\+ b\+ 1)_{n-1}(a\+ 1/2)_n(b\+ 1/2)_n} },
\label{sz-ko} \\
  l_{2n}   & = & 2n{a\m b \over 2n\+ a\+ b}\k_{2n},
\label{sz-le} \\
  l_{2n-1} & = & (2n\m 1){a\m b \over 2n\+ a\+ b\m 1}\k_{2n-1},
\label{sz-lo} \\
 \f_{2n}(0) & = & {a\+ b \over 2n\+ a\+ b}\k_{2n},
\label{sz-phie} \\
 \f_{2n-1}(0) & = & {a\m b \over 2n\+ a\+ b\m 1}\k_{2n-1}.
\label{sz-phio}
\eea
The three term recurrences are then
\begin{multline}
   2(a\m b)\sqrt{n(n\+ a\+ b)}\f_{2n}(z) 
 + 2(a\+ b)\sqrt{(n\+ a\m 1/2)(n\+ b\m 1/2)}z\f_{2n-2}(z) \\
 = [(a\+ b)(2n\+ a\+ b\m 1) + (a\m b)(2n\+ a\+ b)z] \f_{2n-1}(z),
\label{sz-3term1}
\end{multline}
and
\begin{multline}
   2(a\+ b)\sqrt{(n\+ a\m 1/2)(n\+ b\m 1/2)}\f_{2n-1}(z) 
 + 2(a\m b)\sqrt{(n\m 1)(n\+ a\+ b\m 1)}z\f_{2n-3}(z) \\
 = [(a\m b)(2n\+ a\+ b\m 2) + (a\+ b)(2n\+ a\+ b\m 1)z] \f_{2n-2}(z),
\label{sz-3term2}
\end{multline}
when $ a\neq b $ and both these degenerate to
$ \f_{2n-1}(z) = z\f_{2n-2}(z) $ when $ a=b $.
For reference the reciprocal polynomials are
\bea
 z^{-n} \f^{*}_{2n}(z) = \SA P^{(a-1/2,b-1/2)}_{n}(\dfrac{1}{2}[z\+z^{-1}])
       -\dfrac{1}{2} \SB [z\m z^{-1}]P^{(a+1/2,b+1/2)}_{n-1}(\dfrac{1}{2}[z\+z^{-1}]),
\label{sz-conj1} \\
 z^{1-n} \f^{*}_{2n-1}(z) = \SC zP^{(a-1/2,b-1/2)}_{n}(\dfrac{1}{2}[z\+z^{-1}])
       -\dfrac{1}{2} \SD [z\m z^{-1}]zP^{(a+1/2,b+1/2)}_{n-1}(\dfrac{1}{2}[z\+z^{-1}]).
\label{sz-conj2}
\eea
Using the differential and recurrence relations for the Jacobi polynomials
directly one can find the appropriate coefficient functions for the Szeg\"o
polynomials to be
\bea
 A_{2n-1}(z) & = &
   2\sqrt{(n\+ a\m 1/2)(n\+ b\m 1/2)}
   { a\m b + (a\+ b)z \over (a\m b)(1\m z^2) },
\label{sz-Ao} \\
 B_{2n-1}(z) & = &
   { 4ab + (2n\m 1)\left[a\+ b + (a\m b)z\right]
     \over (a\m b)(1\m z^2) },
\label{sz-Bo} \\
 A_{2n}(z) & = &
   2\sqrt{n(n\+ a\+ b)}
   { a\+ b + (a\m b)z \over (a\+ b)(1\m z^2) },
\label{sz-Ae} \\
 B_{2n}(z) & = &
   2n{ a\m b + (a\+ b)z \over (a\+ b)(1\m z^2) },
\label{sz-Be}
\eea
again when $ a \neq b $.

\noindent{\bf Example 3}. 
Consider the weight function
\bea
w(z) = {1 \over 2\pi I_0(t)}\exp( \dfrac{1}{2}t[z+z^{-1}]), 
\label{mB-w}
\eea
where $I_\nu$ is a modified Bessel function. This system of orthogonal
polynomials has arisen from studies of the length of longest increasing 
subsequences of random words \cite{Ba:De:Jo} and matrix models 
\cite{Pe:Sh}, \cite{Hs}.
In terms of the leading coefficient one has the Toeplitz determinant form
\bea
    \k^2_n(t) = I_0(t){ \det( I_{j-k}(t) )_{0 \leq j,k \leq n-1} \over 
                        \det( I_{j-k}(t) )_{0 \leq j,k \leq n} }.
\label{mB-toeplitz}
\eea
The first few members of this sequence are
\bea
  \k^2_1 & = & { I^2_0(t) \over I^2_0(t)-I^2_1(t) },
  \\
  {\f_1(0) \over \k_1} & = & -{I_1(t) \over I_0(t)},
  \\
  \k^2_2 & = & { I_0(t)(I^2_0(t)-I^2_1(t)) \over
   (I_0(t)-I_2(t))\left[I_0(t)(I_0(t)+I_2(t))-2I^2_1(t)\right] },
  \\
  {\f_2(0) \over \k_2} & = &
  { I_0(t)I_2(t)-I^2_1(t) \over I^2_1(t)-I^2_0(t) }.
\label{mB-Tp123}
\eea
Gessel \cite{Ge} has found the exact power series expansions in $ t $
for the first three determinants which appear in the above coefficients.
Some recurrence relations for the corresponding coefficients of the monic
version of these orthogonal polynomials have been known \cite{Pe:Sh},
\cite{Hs}, \cite{Tr:Wi} and we derive the equivalent results for $ \k_n $, etc. 
\begin{lem}[\cite{Pe:Sh}]
The reflection coefficient $ r_n(t) \equiv \f_n(0)/\k_n $ for the modified 
Bessel orthogonal polynomial system satisfies a form of the discrete 
Painl\'eve II equation, namely the recurrence relation
\bea
  -2{n \over t}{r_n \over 1-r^2_n} = r_{n+1}+r_{n-1},
\label{mB-RR}
\eea
for $ n \geq 1 $ and $ r_0(t) = 1 $, $ r_1(t) = -I_1(t)/I_0(t) $.
\end{lem}
{\bf Proof.} Firstly we make a slight redefinition of the external field 
$ w(z) = \exp(-v(z+1/z)) $ for convenience. 
Employing integration by parts we evaluate
\bea
  & &-\int v'(\z+1/\z)(1-1/\z^2)
   \f_{n+1}(\z)\overline{\f_n(\z)} w(\z) \frac{d\z}{i\z} \nonumber\\
  &=&
   \int \left[ \f_{n+1}(\z)\overline{\z}^2\overline{\f_{n}^{\pri}(\z)}
              +\f_{n+1}(\z)\overline{\z}\overline{\f_{n}(\z)}
              -\f_{n+1}^{\pri}(\z)\overline{\f_{n}(\z)} \right]
   w(\z){d\z \over i\z} \nonumber\\
  &=&
  (n\+ 1)\left[ {\k_n \over \k_{n+1}}-{\k_{n+1} \over \k_{n}} \right],
\label{lhs-mB-RR}
\eea
for general external fields $ v(z) $
using (\ref{uc-ortho}) and (\ref{uc-coeff}) in a similar way to the proof of 
Theorem 2.1. However in this case $ v'(\z+1/\z) = - t/2 $, a direct 
evaluation of the left-hand side yields
\bea
  -\frac{1}{2} t\left( {l_n \over \k_{n+1}}-{\k_n l_{n+2} \over 
\k_{n+1}\k_{n+2}}
          \right),
\label{rhs-mB-RR}
\eea
and simplification of this equality in terms of the defined ratio and use of
(\ref{uc-l}) gives the above result.

There is also a differential relation satisfied by these coefficient functions
or equivalently a differential relation in $ t $ for the orthogonal polynomials
themselves \cite{Hs}, \cite{Tr:Wi}.
\begin{lem}
The modified Bessel orthogonal polynomials satisfy the differential relation
\bea
 2{d \over dt}\f_n(z) =
     \left[ {I_1(t) \over I_0(t)}+{\f_{n+1}(0) \over \k_{n+1}}
                                  {\k_n \over \f_n(0)}
     \right] \f_n(z)
    -{\k_{n-1} \over \k_n} \left[ 1+{\f_{n+1}(0) \over \k_{n+1}}
                                  {\k_n \over \f_n(0)}z
                           \right] \f_{n-1}(z),
\label{mB-DR}
\eea
for $ n \geq 1 $ and $ {d \over dt}\f_0(z) = 0 $. The differential equations
for the coefficients are
\bea
  {2 \over \k_n}{d \over dt}\k_n &=&
  {I_1(t) \over I_0(t)}
  +{\f_{n+1}(0) \over \k_{n+1}}{\f_{n}(0) \over \k_{n}},
\label{mB-coeff-DR1} \\
  {2 \over \f_n(0)}{d \over dt}\f_n(0) &=&
  {I_1(t) \over I_0(t)}
 +{\f_{n+1}(0) \over \k_{n+1}}{\k_{n} \over \f_{n}(0)}
 -{\f_{n-1}(0) \over \f_n(0)}{\k_{n-1} \over \k_n},
\label{mB-coeff-DR2}
\eea 
for $ n \geq 1 $.
\end{lem}
{\bf Proof.} Differentiating the orthonormality relation (\ref{uc-ortho}) with 
respect to $ t $ one finds from the orthogonality principle for $ m \leq n-2 $
that
\bea
    {d \over dt}\f_n(z) + \frac{1}{2} z\f_n(z) = 
    a_n\f_{n+1}(z)+b_n\f_n(z)+c_n\f_{n-1}(z)
\eea
for some coefficients $ a_n, b_n, c_n $. The first coefficient is immediately
found to be $ a_n = \frac{1}{2} \k_n/\k_{n+1} $. Consideration of the 
differentiated
orthonormality relation for $ m = n-1 $ sets another coefficient,
$ c_n = - \frac{1}{2} \k_{n-1}/\k_n $, while the case of $ m = n $ leads to
$ b_n = \frac{1}{2} I_1(t)/I_0(t) $. Finally use of the three-term recurrence
(\ref{uc-3term}) allows one to eliminate $ \f_{n+1}(z) $ in favor of 
$ \f_n(z), \f_{n-1}(z) $ and one arrives at (\ref{mB-DR}). The differential
equations for the coefficients $ \k_n, \f_n(0) $ in 
(\ref{mB-coeff-DR1},\ref{mB-coeff-DR2})
follow from reading off the appropriate terms of (\ref{mB-DR}).

Use of the recurrence relation and the differential relations will allow us
to find a differential equation for the coefficients, and thus another
characterization of the coefficients.
\begin{lem}
The reflection coefficient $ r_n(t) $ satisfies the following second order
differential equation
\bea
   {d^2 \over dt^2}r_n =
        \frac{1}{2}\left( {1 \over r_n+1}+{1 \over r_n-1} \right)
        \left({d \over dt}r_n\right)^2 
      - {1 \over t}{d \over dt}r_n - r_n(1-r^2_n)
      + {n^2 \over t^2}{r_n \over 1-r^2_n},
\label{mB-2DE}
\eea
with the boundary conditions determined by the expansion
\bea
  r_n(t) \underset{t \to 0}{\sim}
  {\left(- \dfrac{1}{2}t\right)^n \over n!}\left\{
    1- {1 \over n\+ 1}\frac{1}{4} t^2
                     + {\rm O}(t^4) \right\},
\label{mB-2DE-bc}
\eea
for $ n \geq 1 $. The coefficient $ r_n $ is related by 
\bea
    r_n(t) = {z_n(t)+1 \over z_n(t)-1},
\label{mB-xfm}
\eea
to $ z_n(t) $ which satisfies the Painl\'eve transcendent P-V equation with
the parameters
\bea
   \alpha = -\beta = {n^2 \over 8}, \quad \gamma = 0, \quad \delta = -2.
\label{mB-PV}
\eea
\end{lem}
{\bf Proof.} Subtracting the relations (\ref{mB-coeff-DR1},\ref{mB-coeff-DR2})
leads to the simplified expression
\bea
  r_{n+1}-r_{n-1} = {2 \over 1-r^2_n}{d \over dt}r_n,
\label{mB-RR1}
\eea
which should be compared to the recurrence relation, in a similar form
\bea
  r_{n+1}+r_{n-1} = -2{n \over t}{r_n \over 1-r^2_n}.
\label{mB-RR2}
\eea
The differential equation (\ref{mB-2DE}) is found by combining these latter
two equations and the identification with the P-V can be easily verified.

As a consequence of the above we find that the coefficients for the modified
Bessel orthogonal polynomials can be determined by the Toeplitz determinant
(\ref{mB-toeplitz}), by the recurrence relations (\ref{mB-RR2}) or by the 
differential equation (\ref{mB-2DE}). An example of the use of this last method
we note
\bea
   \k^2_n(t) = I_0(t)\left[ 1-r^2_n(t) \right]^{-1/2}
               \exp\left(-n\int^t_0 {ds \over s}{r^2_n(s) \over 1-r^2_n(s)}
                   \right).
\label{mB-quad}
\eea

We now indicate how to find the coefficients of the differential relations,
$ A_n(z), B_n(z) $ and observe that
\bea
\cq =  -\frac{t}{2} \left[\frac{1}{z\z^2} + \frac{1}{z^2\z} \right]. \nonumber 
\label{mB-logw}
\eea
The above relationship and (\ref{uc-diff3}) yield
\bea
\f_n^{\pri}(z) & = &
     \frac{\k_{n-1}}{\k_n}\f_{n-1}(z) \left[n+ {t \over 2z} \right]
  +  {t \over 2}\frac{\k_{n-1}}{\f_n(0)}\f_{n-1}(z)
     \int_{|\z| =1} \f_n(\z)\overline{\z \f_n^*(\z)} w(\z) \frac{d\z}{i\z} \\
  &{}& \quad
  + {t \over 2z}\f_n(z) \int_{|\z| =1} \f_n(\z)
    \overline{\z} \,  \left[\overline{\f_n(\z)} - \frac{\k_n}{\f_n(0)}
    \overline{\f_n^*(\z)}\right] w(\z) \frac{d\z}{i\z}.  \nonumber 
\label{mB-diffr}
\eea
Easy calculations using (\ref{uc-coeff}) give 
\bea
   \z \left[ \f_n(\z) - \frac{\k_n}{\overline{\f_n(0)}} \f_n^*(\z) \right]
   = - \frac{\k_{n-1}}{\k_n}
       \frac{\overline{\f_{n-1}(0)}}{\overline{\f_n(0)}}\f_n(\z) 
     + {\rm lower \; order \; terms}. 
\nonumber
\label{mB-int1}
\eea
and 
\bea
   \z \frac{\f_n^*(\z)}{\overline{\f_n(0)}}
   = \frac{\f_{n+1}(\z)}{\k_{n+1}} 
     + \left[\frac{\k_nl_n-\k_{n-1}l_{n-1}}{\k_n|\f_n(0)|^2}
             - \frac{l_{n+1}}{k_{n+1}\k_n}\right] \f_n(\z) 
     + {\rm lower \; order \; terms}.\nonumber 
\label{mB-int2}
\eea
These identities together with (\ref{mB-diffr}) establish the 
differential-difference relation 
\bea
   \f_n^{\pri}(z) &=& 
   \frac{\k_{n-1}}{\k_n}
   \left[ n + \frac{t}{2z}
          + \frac{t}{2}\frac{\k_{n-1}}{\k_n}\frac{\f_{n-1}(0)}{\f_n(0)}
          - \frac{t}{2}\frac{\overline{\f_{n+1}(0)}\f_n(0)}{\k_{n+1}\k_n}
   \right] \f_{n-1}(z) \\
   &{}& \quad
     - \frac{t}{2z}\frac{\k_{n-1}}{\k_{n}}\frac{\f_{n-1}(0)}{\f_n(0)}
       \f_n(z)\ . \nonumber 
\label{mB-diff}
\eea

\bigskip

\setcounter{section}{3}
\setcounter{thm}{0}
\setcounter{equation}{0}

\noindent{\bf 3. Functional Equation and Zeros}.
In this section we continue the development of the previous discussion of the
differential relations satisfied by orthogonal polynomials to find a
functional equation and its relationship to the zeros of the polynomials.
Expressing the second order differential equation (\ref{uc-2ode1}) in terms of 
the coefficient functions $ A_n(z) $ and $ B_n(z) $ we have
\begin{multline}
 \f_n''
 + \left\{
    B_n + B_{n-1} - A_n'/A_n
    - {\k_{n-1} \over \k_{n-2}}{A_{n-1} \over z}
    - {\k_n \over \k_{n-2}}{\f_{n-1}(0) \over \f_n(0)}A_{n-1}
   \right\} \f_n' \\
 + \left\{
    B_n' - B_nA_n'/A_n + B_nB_{n-1}
    - {\k_{n-1} \over \k_{n-2}}{A_{n-1}B_n \over z}
   \right.\qquad\qquad \\
   \qquad\qquad\left.
    - {\k_n \over \k_{n-2}}{\f_{n-1}(0) \over \f_n(0)}
      A_{n-1}B_n
    + {\k_{n-1} \over \k_{n-2}}{\f_{n-1}(0) \over \f_n(0)}
      {A_{n-1}A_n \over z}
   \right\} \f_n
 = 0.
\label{uc-2ode1full}
\end{multline}
Now by analogy with the orthogonal polynomials defined on the real line
the coefficient of the $ \f_n' $ term above can be simplified.
\begin{thm}
Given that $ v(z) $ is an meromorphic function in the unit disk then the
following functional equation holds
\bea
    B_n + B_{n-1} - {\k_{n-1} \over \k_{n-2}}{A_{n-1} \over z}
    - {\k_{n} \over \k_{n-2}}{\f_{n-1}(0) \over \f_n(0)}A_{n-1}
  = -(n-1)z^{-1} - v'(z).
\label{uc-AB}
\eea
\end{thm}
{\bf Proof.}
From the definitions (\ref{uc-A},\ref{uc-B}) we start with the following
expression
\begin{multline*}
    B_n + B_{n-1} - {\k_{n-1} \over \k_{n-2}}{A_{n-1} \over z}
    - {\k_{n} \over \k_{n-2}}{\f_{n-1}(0) \over \f_n(0)}A_{n-1} \\
  = -(n-1)\left[\frac{1}{z} + \frac{\k_n}{\k_{n-1}}\frac{\f_{n-1}(0)}{\f_n(0)}
          \right] \\
  \qquad + i\int \cq\left\{ -\f_n\overline{\f_n}
                          +\frac{\k_n}{\f_n(0)}\f_n\overline{\f^*_n}
                          -\f_{n-1}\overline{\f_{n-1}}
                          -\frac{\k_n}{\f_n(0)}\z\f_{n-1}\overline{\f^*_{n-1}}
                  \right\}w(\z)d\z \\
  - i\frac{\k_n}{\f_n(0)}
     \int [v'(z)-v'(\z)]\f_{n-1}\overline{\f^*_{n-1}}w(\z)d\z.
\end{multline*}
Employing the recurrences (\ref{uc-rec2},\ref{uc-rec1}), and the relation
amongst coefficients (\ref{uc-k}) one can show that the factor in the 
first integral on the right-hand side above is
\bea
   -\f_n\overline{\f_n}
   +\frac{\k_n}{\f_n(0)}\f_n\overline{\f^*_n}
   -\f_{n-1}\overline{\f_{n-1}}
   -\frac{\k_n}{\f_n(0)}\z\f_{n-1}\overline{\f^*_{n-1}}
  = -\f_n\overline{\f_n}+\f^*_n\overline{\f^*_n}.
\nonumber
\eea
Now since $|\z|^2 = 1$, one can show that the right-hand side of the above
is zero from the Christoffel-Darboux sum (\ref{uc-cd}). Consequently our
right-hand side is now
\bea
&{}& - (n-1)\left[\frac{1}{z} + \frac{\k_n}{\k_{n-1}}\frac{\f_{n-1}(0)}{\f_n(0)}
        \right]\nonumber\\
 &{}& \quad  - i\frac{\k_n}{\f_n(0)}\left[
           v'(z)\int \f_{n-1}\overline{\f^*_{n-1}}w(\z)d\z
           - \int v'(\z)\f_{n-1}\overline{\f^*_{n-1}}w(\z)d\z \right].
\nonumber
\eea
Taking the first integral in this expression and using the recurrence
(\ref{uc-rec2}) and the decomposition 
$ \z\f_{n-1} = \k_{n-1}/\k_n\f_n + \pi_{n-1} $ where $\pi_{n} \in \Pi_{n} $, 
$\Pi_n$ being  the space
of polynomials of degree at most $n$, we find it reduces to
$ -i\f_n(0)/\k_n $ from the normality of the orthogonal polynomials.
Considering now the second integral above we integrate by parts and are left
with
\bea
  \int  \f'_{n-1}\overline{\f^*_{n-1}}w(\z)d\z
  + \int \f_{n-1}\overline{\f^*_{n-1}}'w(\z)d\z,
\nonumber
\eea
and the first term here must vanish as $ \f^*_{n-1} $ can be expressed in terms
of $ \f_{n-1}, \f_n $ from (\ref{uc-rec2}) but $ \f'_{n-1} \in \Pi_{n-2} $.
The remaining integral, the second one above, can be treated in the following
way. Firstly express the conjugate polynomial in terms of the polynomial
itself via (\ref{uc-rec1}) and employ the relation for its derivative
(\ref{uc-pf1}). Further noting that 
$ \z\f'_{n-1} = (n-1)\f_{n-1}+\pi_{n-2} $,
$ \z\f_{n-2} = \k_{n-2}/\k_{n-1}\f_{n-1}+\pi_{n-2} $, and
$ \z^2\f'_{n-2} = (n-2)\k_{n-2}/\k_{n-1}\f_{n-1}+\pi_{n-2} $ along with
the orthonormality relation, the final integral is nothing but 
$ -i(n-1)\f_{n-1}(0)/\k_{n-1} $.
Combining all this the final result is (\ref{uc-AB}).

\noindent{\bf Remark 1.}
The zeros of the orthogonal polynomial $ \f_n(z) $ are denoted by 
$ \{z_j\}_{1\leq j \leq n} $, and are confined to the convex Hull of the 
support of the measure, namely to be strictly confined within the unit circle
$ |z| < 1 $. One can construct a real function $ |T(z_1, \ldots ,z_n)| $ from
\bea
  T(z_1, \ldots ,z_n) =
  \prod^{n}_{j=1} z^{-n+1}_{j} {e^{-v(z_j)} \over A_{n}(z_j)}
  \prod_{1 \leq j < k \leq n} (z_j-z_k)^2,
\label{uc-T}
\eea
such that the zeros are given by the stationary points of this function. One
might also interpret this function as a total energy function for
$ n $ mobile unit charges in the
unit disk interacting with a one-body confining potential,
$ v(z) + \ln A_n(z) $, an attractive logarithmic potential with a charge
$ n-1 $ at the origin, $ (n-1)\ln z $, and repulsive logarithmic two-body 
potentials, $ -\ln (z_i-z_j) $, between pairs of charges. 
However all the stationary points are saddle-points, a natural consequence of
analyticity in the unit disk.
That such this function exhibits stationary properties at the zeros can be seen
by considering the second order differential equation which in view of the
above theorem has the coefficient $ P(z) $, namely
\bea
  P(z) = -(n-1)z^{-1} - v'(z) - A_n'/A_n.
\label{uc-P}
\eea
This function is a perfect differential and consequently the one-body 
potential can be constructed from its integral, via the Stieltjes 
argument. Or alternatively one can show that the conditions for the
stationary points of function $ T(z_1, \ldots ,z_n) $ above lead to a system of
equations
\bea
   -v'(z_j) - {A_n'(z_j) \over A_n(z_j)} - {n-1 \over z_j}
   + 2 \sum_{1 \leq k \leq n, k \neq j}{1 \over z_j-z_k} = 0
   \quad j = 1, \ldots ,n.
\label{uc-min}
\eea
Then the pairwise sum can be represented in terms of the polynomial
$ f(z) = \prod^{n}_{j=1}(z-z_j) $ thus
\bea
   2 \sum_{1 \leq k \leq n, k \neq j}{1 \over z_j-z_k}
   = {f''(z_j) \over f'(z_j)},
\label{uc-2body}
\eea
and we have the $ n $ conditions expressed as
\bea
   f''(z_j) + \left\{ -{n-1 \over z_j} - v'(z_j) - {A_n'(z_j) \over A_n(z_j)}
              \right\}f'(z_j) = 0 ,\quad \forall j = 1, \ldots ,n.
\label{uc-2odeT}
\eea
The result then follows.

\noindent{\bf Remark 2.}
The functional equation (\ref{uc-AB}) actually implies a very general
recurrence relation on the orthogonal system coefficients 
$ \k_n, \f_n(0) $. In general if it is possible to relate the differential
recurrence coefficients $ A_n, B_n $ to these polynomial coefficients, then
the functional equation dictates that equality holds for all $ z $, and
thus for independent terms in $ z $. For rational functions this can be 
applied to the coefficients of monomials in $ z $.

\noindent{\bf Remark 3.}
There is another way of deriving the functional equation (\ref{uc-AB}) which
we now describe. Equation (\ref{uc-2ode1}) is one way of expressing the second
order differential equation for the orthogonal polynomials, however one can
perform the elimination in the opposite order and find
\bea
L_{n+1,1}\left(\frac{z}{A_n(z)}L_{n+1,2}\right) \f_n(z) =
 \frac{\k_{n}\f_{n}(0)}{\k_{n-1}\f_{n+1}(0)} A_{n+1}(z)\, \f_n(z), 
\label{uc-2ode2}
\eea 
which written out in full is
\begin{multline}
 \f_n''
 + \left\{
    B_{n+1} + B_{n} - A_n'/A_n
    - {\k_{n} \over \k_{n-1}}{A_{n} \over z}
    - {\k_{n+1} \over \k_{n-1}}{\f_{n}(0) \over \f_{n+1}(0)}A_{n}
   + {1 \over z} \right\} \f_n' \\
 + \left\{
    B_n' - B_nA_n'/A_n + B_{n+1}B_{n}
    - {\k_{n} \over \k_{n-1}}{A_{n}B_{n+1} \over z}
   \right. \\
   \quad\left.
    - {\k_{n+1} \over \k_{n-1}}{\f_{n}(0) \over \f_{n+1}(0)}
      A_{n}B_{n+1}
    + {\k_{n} \over \k_{n-1}}{\f_{n}(0) \over \f_{n+1}(0)}
      {A_{n}A_{n+1} \over z} + {B_n \over z} 
    - {\k_{n+1} \over \k_{n-1}}{\f_{n}(0) \over \f_{n+1}(0)}
      {A_n \over z}
   \right\} \f_n
 = 0.
\label{uc-2ode2full}
\end{multline}
Given that the coefficient $ P(z) $ for these two forms 
(\ref{uc-2ode1full},\ref{uc-2ode2full}) must be identical we have an
inhomogeneous first order difference equation, whose solution is
\bea
    B_n + B_{n-1} - {\k_{n-1} \over \k_{n-2}}{A_{n-1} \over z}
    - {\k_{n} \over \k_{n-2}}{\f_{n-1}(0) \over \f_n(0)}A_{n-1}
  = -(n-1)z^{-1} + \text{function of $ z $ only}.
\label{uc-ABnew}
\eea
This function can be simply evaluated by setting $ n = 1 $, evaluating
the integrals after noting $ B_0 = 0 $ and the cancellations, and 
yields the result $ -v'(z) $.

\noindent{\bf Example 1}.
We can verify that the general form for the $ T $-function is correct in
the case of the circular Jacobi polynomials by a direct evaluation
\bea
   T(z_1,\ldots ,z_n) =
   \prod_{j=1}^{n} z^{1-n-a}_j (1-z_j)^{a+1} (z_j-1)^{a}
   \prod_{1 \leq j < k \leq n} (z_j-z_k)^2,
\label{cj-T}
\eea
where we have used the identity
\bea
   |1-z|^{2a} = (1-z)^a(1-1/z)^a = z^{-a}(1-z)^a(z-1)^a,
\label{cj-map}
\eea
on $ |z|=1 $ to suitably construct a locally analytic weight function.
One can show that the stationary points for this problem are the solution to 
the set of equations
\bea
   {1-n-a \over z_j} -{2a+1 \over 1-z_j}
   + 2\sum_{j \neq k} {1 \over z_j-z_k} = 0 , \quad 1 \leq j \leq n,
\label{cj-min}
\eea
so that the polynomial $ f(z) = \prod_{j=1}^n (z-z_j) $
satisfies the relations
\bea
    f''(z_j) + f'(z_j)\left\{
    {1\m n\m a \over z_j} - {2a\+ 1 \over 1-z_j} \right\} = 0.
\label{cj-2ode2}
\eea
Consequently we find that
\bea
   z(1-z) f''(z) + f'(z)\left\{
    (1\m n\m a)(1-z) - (2a\+ 1)z  \right\} + Qf(z) = 0,
\label{cj-2ode3}
\eea
for some constant $ Q $ independent of $ z $, but possibly dependent on
$ n $ and $ a $ and is identical to the second order ODE (\ref{cj-2ode}).

\noindent{\bf Example 2}.
Using the expressions (\ref{sz-Ao}-\ref{sz-Be}) one can verify that the 
identity (\ref{uc-AB}) holds and in particular becomes
\bea
  B_{n}+B_{n-1}-{\k_{n-1} \over \k_{n-2}}{A_{n-1} \over z}
  -{\k_{n} \over \k_{n-2}}{\f_{n-1}(0) \over \f_{n}(0)}A_{n-1}
  = -{n\m 1 \over z} - {a\+ b \over z} - {2a \over 1\m z} + {2b \over 1\+ z},
\label{sz-AB}
\eea
for both the odd and even sequences. Consequently the coefficients in the
second order differential equation are
\bea
  P_{n}(z) & = & -{n\+ a\+ b\m 1 \over z} - {2a\+ 1 \over 1\m z}
             + {2b\+ 1 \over 1\+ z} - {a\pm b \over a\mp b + (a\pm b)z},
\label{sz-P} \\
  Q_{2n}(z) & = & 2n
  { a(a\+ 1)(1\+ z)^2 - b(b\+ 1)(1\m z)^2 \over 
    z(1\m z^2)[a\+ b + (a\m b)z] },
\label{sz-Qe} \\
  Q_{2n-1}(z) & = & 
  { (2n-1)[a(a\+ 1)(1\+ z)^2 + b(b\+ 1)(1\m z)^2 -2ab(1\m z^2)] + 4ab \over 
    z(1\m z^2)[a\m b + (a\+ b)z] }.
\label{sz-Qo}
\eea
 
Similarly we can verify that the general form for the $ T $-function is correct
in the case of the Szeg\"o polynomials by using the identity
\bea
   |1-z|^{2a}|1+z|^{2b} = z^{-a-b}(1-z)^a(z-1)^a(1+z)^{2b},
\label{sz-map}
\eea
to suitably analytically continue the weight function.
The stationary points for this problem are the solution to the following set 
of equations
\bea
   {1\m n\m a\m b \over z_j} -{2a\+1 \over 1-z_j} +{2b\+1 \over 1+z_j}
   -{a\mp b \over a\pm b + (a\mp b)z_j}
   + 2\sum_{j \neq k} {1 \over z_j-z_k} = 0 , \quad 1 \leq j \leq n,
\label{sz-min}
\eea
such that the polynomial $ f(z) = \prod_{j=1}^n (z-z_j)$
satisfies the relations
\bea
    f''(z_j) + f'(z_j)\left\{
    {1\m n\m a\m b \over z_j} - {2a\+ 1 \over 1-z_j} 
   +{2b\+1 \over 1+z_j} -{a\mp b \over a\pm b + (a\mp b)z_j} \right\} = 0.
\label{sz-2ode2}
\eea
Finally we find that
\bea
   f''(z) + f'(z)\left\{ 
    {1\m n\m a\m b \over z} - {2a\+ 1 \over 1-z} 
   +{2b\+1 \over 1+z} -{a\mp b \over a\pm b + (a\mp b)z} \right\}
   + Q(z)f(z) = 0,
\label{sz-2ode3}
\eea
for some constant $Q$ independent of $z$, but possibly dependent on
$ n $ and $ a $. The coefficient of the first derivative term is identical to 
the expression for $ P(z) $ in (\ref{sz-P}).

\noindent{\bf Example 3}.
One can also verify the functional relation (\ref{uc-AB}) for the modified 
Bessel orthogonal polynomials. Forming the left-hand side of this identity
we find this reduces to
\begin{multline}
  B_n + B_{n-1} - {\k_{n-1} \over \k_{n-2}}{A_{n-1} \over z}
  - {\k_n \over \k_{n-2}}{\f_{n-1}(0) \over \f_{n}(0)} A_{n-1} \\
  = - {n-1 \over z} - {t \over 2z^2}
    - (n-1){\k_n \over \k_{n-1}}{\f_{n-1}(0) \over \f_{n}(0)}
    - {t \over 2}{\k_n\k_{n-2} \over \k^2_{n-1}}{\f_{n-2}(0) \over \f_{n}(0)}
    + {t \over 2}{\f^2_{n-1}(0) \over \k^2_{n-1}}.
\label{}
\end{multline}
Now the last three terms on the right-hand side of the above equation
simplify to $ t/2 $ using the recurrence relation (\ref{mB-RR}), showing that
the general functional relation holds. In fact, as remarked earlier, this
relation implies the recurrence relation itself.

\bigskip

\setcounter{section}{4}
\setcounter{thm}{0}
\setcounter{equation}{0}

\noindent{\bf 4. $q$-Difference Equations}. 
Our first result is a $q$-analogue of Theorem 2.1. 
\begin{thm}
If $w(z)$ is analytic in the ring  $q < |z| < 1$ and is continuous 
on its boundary then  
\bea
(D_q\f_n)(z) &= & \frac{\k_{n-1}}{\k_n} \frac{1\m q^n}{1\m q} 
\f_{n-1}(z)  
 -i\f_n^*(z) \int_{|\z|=1} \frac{u(\z) - u(qz)}{\z-qz}
  \f_n(\z) 
 \overline{\f_n^*(q\z)}w(\z) d\z \\
 &{}& \quad +i \f_n(z) \int_{|\z|=1} \frac{u(\z) - u(qz)}{\z-qz}
  \f_n(\z)     
 \overline{\f_n(q\z)} w(\z)  d\z.  \nonumber   
\label{qop-diffr}
\eea 
\end{thm}
{\bf Proof}. Expand $D_q\f_n(z)$ in a series of the $\f_n$'s. We get 
\bea
(1\m q) (D_q\f_n)(z) = \int_{|\z|=1} 
\sum_{k=0}^{n-1}\f_k(z)\overline{\f_k(\z)}[\f_n(\z)-\f_n(q\z)] w(\z)
\frac{d\z}{i\z^2}. \nonumber
\label{qop-diff1}
\eea 
Break the above integral as a difference of two integrals involving  
$\f_n(\z)$ and  $\f_n(q\z)$, then in the second integral replace $\z$ 
by $\z/q$. Under such transformation $\overline{\f_k(\z)}$ is transformed to 
$\overline{\f_k(q\z)}$, since $|\z| = 1$.  Furthermore (\ref{qop-logw}) gives 
\bea
w(\z/q) = [1+ \z u(\z)(1-1/q)]w(\z).    
\label{qop-w2}
\eea
Therefore 
\bea
(1\m q) (D_q\f_n)(z) &=&  \int_{|\z|=1} 
\sum_{k=0}^{n-1}\f_k(z)\overline{\z\f_k(\z)}\f_n(\z) w(\z)
\frac{d\z}{i\z}   \nonumber \\ 
&{}& \;    + \int_{|\z|=1} \sum_{k=0}^{n-1}\f_k(z) 
\left[-q \overline{\z\f_k(q\z)}  + u(\z)(1\m q)\overline{\f_k(q\z)}\right]
 \f_n(\z) w(\z) \frac{d\z}{i\z}     \nonumber \\
&=& \frac{\k_{n-1}}{\k_n}\f_{n-1}(z) 
  - q^n \frac{\k_{n-1}}{\k_n}\f_{n-1}(z) \nonumber \\
 &{}& \;   + (1\m q)\int_{|\z|=1} \f_n(\z)u(\z) \sum_{k=0}^{n-1}\f_k(z) 
 \overline{\f_k(q\z)} w(\z)\frac{d\z}{i\z}.  \nonumber   
\label{qop-diff2}
\eea 
The result now follows from (\ref{uc-cd}). 

We next substitute for $\f_n^*(z)$ in (\ref{qop-diffr}) from (\ref{uc-rec1}),
if $ \f_n(0) \neq 0 $, and establish
\bea
(D_q\f_n)(z) = A_n(z)\f_{n-1}(z) - B_n(z) \f_n(z),  
\label{qop-diff}
\eea
with 
\bea
A_n(z) &=&  \frac{\k_{n-1}}{\k_n} \frac{1\m q^n}{1\m q}  
 + i\frac{\k_{n-1}}{\f_n(0)} z
 \int_{|\z|=1} \frac{u(\z) - u(qz)}{\z-qz}  \f_n(\z) 
 \overline{\f_n^*(q\z)}w(\z) d\z
\label{qop-A} \\
B_n(z)  &=& \quad  -i\int_{|\z|=1} \frac{u(\z) - u(qz)}{\z-qz}
  \f_n(\z)\left[\overline{\f_n(q\z)} -\frac{\k_n}{\f_n(0)}
\overline{\f_n^*(q\z)}\right]
\, w(\z)  d\z.    
\label{qop-B}
\eea 
The are the  $q$-analogues of (\ref{uc-A}-\ref{uc-diff}). Here again we set 
\bea
L_{n,1} &=&  D_q +  B_n(z),
\label{qop-L1} \\
L_{n,2} &=& - D_q - B_{n-1}(z) + \frac{A_{n-1}(z)\k_{n-1}}{z\k_{n-2}} 
 + \frac{A_{n-1}(z)\k_{n }\f_{n-1}(0)}{\k_{n-2}\f_n(0)}. 
\label{iqop-L2}
\eea
The ladder operations are 
\bea
L_{n,1} \f_n(z) = A_n(z) \f_{n-1}(z), \quad 
L_{n,2} \f_{n-1}(z) = 
\frac{\f_{n-1}(0)\k_{n-1}}{\f_n(0)\k_{n-2}} \frac{A_{n-1}(z)}{z} \f_n(z). 
\label{qop-ac}
\eea
This results in the $q$-difference  equation
\bea
L_{n,2}\left(\frac{1}{A_n(z)}L_{n,1}\right) \f_n(z) =
 \frac{A_{n-1}(z)}{z}\, \frac{\f_{n-1}(0)\k_{n-1}}{\f_n(0)\k_{n-2}}\f_n(z). 
\label{qop-2ode1}
\eea 

There is also a $q$-analogue of the functional equation (\ref{uc-AB}) which
can be found most simply by exploiting the third Remark to Theorem (3.1).
\begin{thm}
If $ u(z) $ is analytic in the annular region $ q < |z| < 1 $ then the following
functional equation for the coefficients $ A_n(z), B_n(z) $ holds 
\begin{multline}
   B_n + B_{n-1} -{\k_{n-1} \over \k_{n-2}}{A_{n-1} \over z}
   - {\k_n \over \k_{n-2}}{\f_{n-1}(0) \over \f_{n}(0)}A_{n-1} \\
   = -{n-1 \over qz}-{u(qz) \over q}
     -{1\m q \over q}\sum_{j=0}^{n-1} \left[
           B_{j+1}-{\k_{j} \over \k_{j-1}}{A_{j} \over z} \right].
\label{qop-FEq}
\end{multline}
\end{thm}
\noindent{\bf Proof.}
Two alternative forms of the second order $q$-difference equation are possible,
namely (\ref{qop-2ode1}) and the following,
\bea
L_{n+1,1}\left({z \over A_n(z)}L_{n+1,2}\right)\f_n(z) 
  = {\k_n\f_n(0) \over \k_{n-1}\f_{n+1}(0)}A_{n+1}(z)\f_n(z).
\label{qop-2ode2}
\eea
These two equations, written out in full are, respectively
\begin{multline}
 D^2_q\f_n(z) \\
 + \left\{
    B_{n}(qz) + {A_n(qz) \over A_n(z)}B_{n-1}(z) - {D_qA_n(z) \over A_n(z)}
   \right. \\
   \qquad\qquad\left.
    - {\k_{n-1} \over \k_{n-2}}{A_n(qz)A_{n-1}(z) \over A_n(z)z}
    - {\k_{n} \over \k_{n-2}}{\f_{n-1}(0) \over \f_{n}(0)}
      {A_n(qz)A_{n-1}(z) \over A_n(z)}
   \right\} D_q\f_n(z) \\
 + \left\{
    D_qB_n(z) - {B_n(z) \over A_n(z)}D_qA_n(z) 
    + {A_n(qz) \over A_n(z)}B_{n}(z)B_{n-1}(z)
    - {\k_{n-1} \over \k_{n-2}}
      {A_n(qz) \over A_n(z)}{A_{n-1}(z)B_{n}(z) \over z}
   \right. \\
   \qquad\left.
    - {\k_{n} \over \k_{n-2}}{\f_{n-1}(0) \over \f_{n}(0)}
      {A_n(qz) \over A_n(z)}A_{n-1}(z)B_{n}(z)
    + {\k_{n-1} \over \k_{n-2}}{\f_{n-1}(0) \over \f_{n}(0)}
      {A_{n}(qz)A_{n-1}(z) \over z}
   \right\} \f_n(z)
 = 0,
\label{qop-2ode1full}
\end{multline}
and
\begin{multline}
 D^2_q\f_n(z) \\
 + \left\{
    {A_n(qz) \over qA_n(z)}B_{n+1}(z) + B_{n}(qz) - {D_qA_n(z) \over qA_n(z)}
   \right. \\
   \qquad\qquad\left .
    - {\k_{n} \over \k_{n-1}}{A_{n}(qz) \over qz}
    - {\k_{n+1} \over \k_{n-1}}{\f_{n}(0) \over \f_{n+1}(0)}A_{n}(qz)
   + {1 \over qz} \right\} D_q\f_n(z) \\
 + \left\{
    D_qB_n(z) - {B_n(z) \over qA_n(z)}D_qA_n(z) 
    + {A_n(qz) \over qA_n(z)}B_{n+1}(z)B_{n}(z)
    - {\k_{n} \over \k_{n-1}}{A_{n}(qz)B_{n+1}(z) \over qz}
   \right.\qquad\qquad \\
   \qquad\qquad
    - {\k_{n+1} \over \k_{n-1}}{\f_{n}(0) \over \f_{n+1}(0)}
      {A_{n}(qz)B_{n+1}(z) \over q}
    + {\k_{n} \over \k_{n-1}}{\f_{n}(0) \over \f_{n+1}(0)}
      {A_{n}(qz)A_{n+1}(z) \over qz} \\
   \left.
    + {B_n(z) \over qz} 
    - {\k_{n+1} \over \k_{n-1}}{\f_{n}(0) \over \f_{n+1}(0)}
      {A_n(qz) \over qz}
   \right\} \f_n(z)
 = 0.
\label{qop-2ode2full}
\end{multline}
A comparison of the coefficients of the first $q$-difference terms leads to
the difference equation
\begin{multline}
   {1 \over q}B_{n+1}(z)-B_{n-1}(z)
   -{\k_n \over \k_{n-1}}{A_n(z) \over qz}
   +{\k_{n-1} \over \k_{n-2}}{A_{n-1}(z) \over z} \\
   -{\k_{n+1} \over \k_{n-1}}{\f_{n}(0) \over \f_{n+1}(0)}A_n(z)
   +{\k_{n} \over \k_{n-2}}{\f_{n-1}(0) \over \f_{n}(0)}A_{n-1}(z)
   = -{1 \over qz}.
\label{qop-FEdiff}
\end{multline}
Using the results for the first coefficients 
\bea
   B_1(z) & = & -u(qz) - {\f_1(qz) \over \f_1(0)}M_1(z), \\
   {A_0(z) \over \k_{-1}} & = & -zM_1(z),
\label{qop-firstAB}
\eea
with
\bea
   M_1(z) \equiv \int_{|\z|=1}
               \z {u(\z)-u(qz) \over \z-qz} w(\z)\frac{d\z}{i\z},
\label{qop-Mdefn}
\eea
this difference equation can be summed to yield the result in (\ref{qop-FEq}).

In the next example we will follow the notation and terminology in \cite{Ga:Ra}
and \cite{An:As:Ro}. The $q$-shifted factorials are 
\bea
(a; q)_0 := 1, \quad (a;q)_n := \prod_{k=1}^n (1-aq^{k-1}), \quad n = 1, 
\cdots {\rm or} \; \infty, \nonumber 
\label{q-fact}
\eea
while the multi-shifted factorials are 
\bea
(a_1, a_2, \cdots, a_m;q)_n := \prod_{k=1}^m(a_k;q)_n. \nonumber 
\label{q-mfact}
\eea

Define the inner product 
\bea
(f, g) = \int_{|z| =1} f(\z)\overline{g(\z)} w(\z) \frac{d\z}{\z}. 
\label{q-ip}
\eea
With respect to this inner product the adjoint of $L_{n,1}$ is 
\bea
(L_{n,1}^*f) (z) = z^2[q- (1\m q)z \overline{u(z)}]  D_q f(z)
+ z f(z) + \left[\overline{B_n(z)} + \overline{u(z)}\right]f(z), 
\label{q-adL1}
\eea
provided that $w(z)$ is analytic in $q < |z| < 1$ and is continuous on 
$|z| =1$ and $|z| =q$.   
The proof follows from  the definition of $D_q$ and the fact 
$\overline{g(\z)} = \overline{g}(1/\z)$, when $|\z| =1$. Observe that as 
$q \to 1$, the right-hand side of (\ref{q-adL1}) tends to the right-hand side of
(\ref{uc-adL1}), as expected.  

\noindent{\bf Example}. 
Consider the Rogers-Szeg\H{o} \cite{Sz1} polynomials $\{{\cal H}_n(z|q)\}$, 
where
\bea
w(z) = \frac{(q^{1/2}z, q^{1/2}/z;q)_\infty}{2\pi (q;q)_\infty}, 
\quad {\cal H}_n(z|q) = \sum_{k=0}^n \frac{(q;q)_nq^{-k/2}z^k}
{(q;q)_k(q;q)_{n-k}}. 
\label{rs-def}
\eea
In this case 
\bea
\f_n(z) = \frac{q^{n/2}}{\sqrt{(q;q)_n}}{\cal H}_n(z|q), 
\quad \f_n(0) = \frac{q^{n/2}}{\sqrt{(q;q)_n}}, 
\quad \k_n = \frac{1}{\sqrt{(q;q)_n}}.
\label{rs-poly}
\eea
It is easy to see that 
\bea
u(z) = \frac{\sqrt{q}}{1\m q} + \frac{qz^{-1}}{1\m q} 
\label{rs-u}
\eea
Thus $[u(\z)-u(qz)]/(\z -qz)$ is $-1/[(1\m q)\z z]$. 
A simple calculation gives 
\bea
(D_q\f_n)(z) = \frac{(1\m q^n)^{3/2}}{1\m q} \f_{n-1}(z) + 
\frac{\k_{n-1}\f_{n-1}(z)}{\f_n(0)(1\m q)} \int_{|\z|=1} \f_n(\z) 
\overline{\f_n^*(q\z)}\, w(\z)\frac{d\z}{i\z}, \nonumber
\label{rs-diffr}
\eea
which simplifies to 
\bea
(D_q\f_n)(z) = \frac{\sqrt{1\m q^n}}{1\m q} \f_{n-1}(z), 
\label{rs-diff}
\eea
since $\k_n \overline{\f_n^*(q\z)}-\f_n(0)\overline{\f_n(q\z)}$ is 
a polynomial of degree $n-1$. 
The functional equation (\ref{rs-diff}) can be verified independently  
by direct computation.  

\bigskip

\setcounter{section}{5}
\setcounter{thm}{0}
\setcounter{equation}{0}

\noindent{\bf 5. Discriminants}. 
Schur \cite{Sc}, \cite[\S 6.71]{Sz2} gave an interesting 
proof of the Stieltjes-Hilbert evaluation of the discriminants of the 
classical orthogonal polynomials of Hermite, Laguerre, and Jacobi. His proof 
relies on a very clever observation. Let $\{p_n(x)\}$ be a
sequence of polynomials satisfying a three term recurrence relation 
\bea
p_{n}(x) =  (a_n x + b_n) p_{n-1}(x) - c_n p_{n-2}(x), \quad n >1, 
\label{op-3term}
\eea
and  the initial conditions
\bea
p_0(x) =1; \quad p_1(x) = a_1 x+ b_1, 
\label{op-init}
\eea
together with the conditions $a_{n-1}c_n \ne 0, n > 1$. Schur \cite{Sc}
observed that  
\bea
\prod_{k=1}^n p_{n-1}(x_{j,n})
= (-1)^{n(n-1)/2}\prod_{k=1}^na_{k}^{n-2k+1}c_k^{k-1}, 
\label{op-schur}
\eea
where $\{x_{j,n}: 1 \le j \le n\}$ is the set of zeros of $p_n(x)$.  

Let $z_{j,n}$ be the zeros of $\f_n(z)$. Following Schur, we let 
\bea
\Delta_n =  \prod_{j=1}^n \f_{n-1}(z_{j,n}). 
\label{op-disc}
\eea
\begin{lem}
The expression $\Delta_n$ is given by
\bea
\Delta_n = \frac{[\f_n(0)]^{n-1}}{\k_n^{n-1}\k_{n-1}^{n}} \, 
\prod_{j=1}^{n-1} \k_j^2, \quad n\ge 2, \quad \Delta_1 = 1.
\label{uc-disc}
\eea
\end{lem}
{\bf Proof}.  It is clear that
\bea
\Delta_n &=& \k_{n-1}^n \prod_{j=1}^n \prod_{k=1}^{n-1}(z_{j,n}- z_{k,n-1}) 
=   \frac{\k_{n-1}^n}{k_n^{n-1}} \prod_{k=1}^{n-1}\f_n(z_{k,n-1})
\label{uc-disc1}
\eea
The recurrence relation (\ref{uc-3term}) we find 
\bea
\f_n(z_{k,n-1})= -\frac{\k_{n-2}\f_n(0)}{\k_{n-1} \f_{n-1}(0)}\; 
z_{k,n-1}\f_{n-2}(z_{k,n-1}).
\label{uc-zero}
\eea
Substituting from (\ref{uc-zero}) into (\ref{uc-disc1}) and applying 
$\f_{n-1}(0) = \k_{n-1}\prod_{j=1}^{n-1}(-z_{j,n-1})$ we establish the two term
recurrence relation   
\bea
\Delta_n = \frac{\k_{n-2}^{n-1} \, [\f_n(0)]^{n-1}}
{\k_{n}^{n-1} \, [\f_{n-1}(0)]^{n-2}}\, \Delta_{n-1}, \quad n > 1. \nonumber
\label{uc-drec}
\eea
By direct computation we find $\Delta_1 = 1$, so the above two term 
recursion implies (\ref{uc-disc}). 

\noindent{\bf Examples}. For the circular Jacobi polynomials the discriminant 
is given by
\bea
\Delta_n = \left( {a \over n\+ a} \right)^{n-1}
      \left[ {(n\m 1)!(2a\+ 1)_{n-1} \over (a\+ 1)^2_{n-1}} \right]^{n/2}
      \prod^{n-1}_{j=1}{(a\+ 1)^2_j \over j!(2a\+ 1)_j},
\label{uc-discCJ}
\eea
while those for the Szeg\"o polynomials are
\bea
\Delta_{2n} &=&
      \left( {a\+ b \over 2n\+ a\+ b} \right)^{2n-1}
      \left[ {(n\m 1)!(a\+ b\+ 1)_{n-1}(a\+ 1/2)_{n}(b\+ 1/2)_{n} 
             \over (a\+ b\+ 1)^2_{2n-1}} \right]^{n}
  \nonumber\\
  & & \times{1 \over (a\+ 1/2)_n(b\+ 1/2)_n}
      \left[ {\prod^{2n-1}_{j=1}(a\+ b\+ 1)_j \over 
              \prod^{n-1}_{l=1}l!(a\+ b\+ 1)_l(a\+ 1/2)_l(b\+ 1/2)_l} \right]^2,
  \\
\Delta_{2n-1} &=&
      \left( {a\m b \over 2n\m 1\+ a\+ b} \right)^{2n-2}
      \left[ {(n\m 1)!(a\+ b\+ 1)_{n-1}(a\+ 1/2)_{n-1}(b\+ 1/2)_{n-1}
             \over (a\+ b\+ 1)^2_{2n-2}} \right]^{n-1/2}
  \nonumber\\
  & & \times(n\m 1)!(a\+ b\+ 1)_{n-1}
      \left[ {\prod^{2n-2}_{j=1}(a\+ b\+ 1)_j \over 
              \prod^{n-1}_{l=1}l!(a\+ b\+ 1)_l(a\+ 1/2)_l(b\+ 1/2)_l} \right]^2.
\label{uc-discSP}
\eea

The resultant of two polynomials $f_n$ and $g_m$ is
\bea
R\{f_n, g_m\} = \gamma^m \prod_{j=1}^n g_m(z_j),
\label{result}
\eea
where $f_n$ is as in (\ref{poly-disc}). Observe that \cite[\S 100]{Di}
\bea
D(f_n) = (-1)^{n(n-1)/2}\gamma^{-1}\, R\{f_n, f_n^{\prime}\}.
\label{res-disc}
\eea
In general let $T$ be a degree reducing operator $T$, that is
 $(Tf)(x)$ is a polynomial of exact 
degree $n-1$ when $f$ has precise degree $n$ and the leading terms in $f$ and 
$Tf$ have the same sign. Define the  generalized discriminant $D(f_n, T)$ by 
\bea
D(f_n, T) := (-1)^{n(n-1)/2} \gamma^{-1} R\{f_n, Tf_n\} = 
(-1)^{n(n-1)/2} \gamma^{n-2} \prod_{j=1}^n (Tf_n)(z_j),
\label{gen-disc}
\eea
for $f_n$ as in (\ref{poly-disc}).  

\begin{thm}
Let $\{\f_n(z)\}$ be orthonormal on the unit circle and assume that $T$ is 
a linear  operator such that 
\bea
T\f_n(z) = A_n(z) \f_{n-1}(z) - B_n(z) \f_n(z).
\label{def-Top}
\eea
Let $\{z_{k,n}:1 \le k \le n\}$ be the zeros of $\f_n(z)$. Then the generalized 
discriminant (\ref{gen-disc}) is given by 
\bea
D(\f_n, T) = \frac{(-1)^{n(n-1)/2}\, [\f_n(0)]^{n-1}}
                  {\k_n\k_{n-1}^{n}} 
\; \prod_{j=1}^{n-1}\k_j^2 \; \prod_{k=1}^n A_n(z_{k,n}), \quad n > 0.
\label{uc-gdisc}
\eea 
\end{thm}
{\bf Proof}. Apply (\ref{def-Top}), (\ref{op-disc}) and (\ref{uc-disc}).

In the case of the orthonormal Rogers-Szeg\H{o} polynomials Theorem 4.2 and, 
(\ref{rs-poly}) and (\ref{rs-diff}) imply the discriminant formula
\bea
D(\f_n, D_q)= D(\f_n,q) = \frac{(-q)^{n(n-1)/2}} {(1\m q)^n}\; (q;q)_n
\prod_{j=1}^{n-1}\left[\frac{1}{(q;q)_j}\right]. 
\label{rs-disc}
\eea
For the Rogers-Szeg\H{o} polynomials we get
\bea
D({\cal H}_n,q) = (-q)^{-n(n-1)/2} \left[{(q;q)_n \over (1\m q)}\right]^{n}
\prod_{j=1}^{n-1}\left[\frac{1}{(q;q)_j}\right]. 
\label{rs-disc2}
\eea
If one is interested in the limiting case $q \to 1$ then we need to rewrite 
(\ref{rs-disc2}) as 
\bea
 D({\cal H}_n,  q) = (-q)^{n(n-1)/2}(1\m q)^{n(n-1)/2}
          \left[{(q;q)_n \over (1\m q)^n}\right]^{n}
          \prod_{j=1}^{n}\left[\frac{(1\m q)^j}{(q;q)_j}\right], 
\label{rs-disc3}
\eea
which shows that $D({\cal H}_n,  q) \to 0$, for $n > 1$,   
when $q \to 1$, as expected since ${\cal H}_n(z|q) \to (1+z)^n$ as $q \to 1$.

\bigskip

{\bf Acknowledgements}. This work started when M. Ismail was visiting the 
University of Melbourne while supported by Omar Foda's research grant. Thanks 
Omar for the hospitality and thanks to Omar and Peter Forrester for many 
enlightening discussions.

\begin{enumerate}
\item
Department of Mathematics, University of South Florida, Tampa, Florida 
33620-5700, U.S.A.
\item
Department of Mathematics and Statistics \& School of Physics, University of 
Melbourne, Victoria 3010, Australia. 
\end{enumerate}


\bigskip  

\begin{thebibliography}{xxx}

\bibitem{An:As:Ro}  G. E. Andrews, R. A. Askey and R. Roy, Special Functions, 
Cambridge University Press, Cambridge, 1999. 
%
\bibitem{Ba} W. Bauldry, {\it Estimates of asymmetric Freud polynomials on 
the real line}, J. Approximation Theory {\bf 63} (1990), 
225--237. 
%
\bibitem{Ba:De:Jo} J. Baik, P. Deift and K. Johansson, {\it On the distribution
of the length of the longest increasing subsequence of random permutations},
J. Amer. Math. Soc. {\bf 12} (1999), 1119-1178.
%
\bibitem{Bo:Cl} S. S. Bonan and D. S. Clark, {\it Estimates of the Hermite 
and the Freud polynomials}, J. Approximation Theory {\bf 63} (1990), 
210--224. 
%
 \bibitem{Ch:Is} Y. Chen and M. E. H. Ismail, {\it Ladder operators and 
 differential equations for orthogonal polynomials},  J. Phys. A  {\bf 30}
 (1997) 7818--7829. 
%
\bibitem{Di} L. E. Dickson, New Course on the Theory of Equations, 
Wiley, New York, 1939.
%
\bibitem{Fo:Ro}  P. J. Forrester and J. B. Rogers, {\it Electrostatics 
and the zeros of the classical orthogonal polynomials}, 
SIAM  J. Math. Anal.  {\bf 17} (1986), 461--468. 
%
\bibitem{Fr}  G. Freud, Orthogonal Polynomials,
Pergamon, Oxford, 1971.
%
\bibitem{Ga:Ra} G. Gasper and M. Rahman, Basic Hypergeometric Series, 
Cambridge University Press, Cambridge, 1990.
%
\bibitem{Ge1}  L. Ya. Geronimus, Orthogonal Polynomials,
Consultants Bureau, New York, 1971.
%
\bibitem{Ge2}  L. Ya. Geronimus, {\it Polynomials orthogonal on a circle and
their applications}, Amer. Math. Soc. Transl., American Mathematical Society, 
Rhode Island, {\bf 3} (1962), 1--78.
%
\bibitem{Ge3}  L. Ya. Geronimus, {\it Orthogonal polynomials}, 
Amer. Math. Soc. Transl., American Mathematical Society, Rhode Island, 
{\bf 108} (1977), 37--130.
%
\bibitem{Ge}  I. M. Gessel, {\it Symmetric Functions and P-Recursiveness},
J. Comp. Theor., Ser. A{\bf 53} (1990), 257--285. 
%
\bibitem{Gr}  F. A. Grunbaum, {\it Variation on a theme of 
Stieltjes and Heine: an electrostatic interpretation of zeros of 
certain polynomials}, J. Comp. Appl. Math. {\bf 99} (1998), 189--194. 
%
\bibitem{Hi} D. Hilbert, {\it \"Uber die Discriminante der in Endlichen 
abbrechenden hypergeometrischen Reihe}, J. f\"ur die reine und angewandte 
Matematik, {\bf 103} (1885), 337--345.
%
\bibitem{Hs} M. Hisakado, {\it Unitary Matrix Models and {P}ainl\'eve III},
Mod. Phys. Lett. A, {\bf 11} (1996), 3001--3010.
%
\bibitem{Ho:Jo} R. A. Horn and C. Johnson, Matrix Analysis, Cambridge 
University Press, Cambridge, 1992. 
%
\bibitem{Is1} M. E. H. Ismail, {\it Discriminants and functions of 
the second kind of orthogonal polynomials}, Results in Mathematics 
{\bf 34} (1998), 132--149. 
%
\bibitem{Is2} M. E. H. Ismail, {\it An electrostatics model for zeros 
of general orthogonal polynomials}, Pacific J. Math. (2000), to appear. 
%
\bibitem{Is3} M. E. H. Ismail, {\it More on electrostatics models for zeros
of orthogonal polynomials}, J. Nonlinear Functional Analysis and Optimization, 
to appear. 
%
\bibitem{Is4} M. E. H. Ismail, {\it Difference equations and  
quantized discriminants for  $q$-orthogonal polynomials}, Advances in 
Applied Math. (2000), to appear. 
%
\bibitem{Ko} T. H. Koornwinder, {\it Orthogonal polynomials with weight 
function}  $(1-x)^\a (1+x)^\b + M\delta(x+1) + N\delta(x-1)$, 
Canadian Math Bull. {\bf 27}(1984), 205--214. 
%
\bibitem{Ma} A. Magnus, MAPA3072A Special topics in approximation theory 
1999-2000: Semi-classical orthogonal polynomials on the unit circle,
{\tt http://www.math.ucl.ac.be/\~{}magnus/}.
%
\bibitem{Ne} P. Nevai, {\it G\'eza Freud, orthogonal polynomials, and 
Christoffel functions. A case study}, J. Approx. Theory {\bf 48} (1986), 3--167.
%
\bibitem{Pa} P.~I.~Pastro, {\it Orthogonal polynomials and some $q$-beta 
integrals of Ramanujan}, J. Math. Anal. Appl. {\bf 112} (1985), 517--540.
%
\bibitem{Pe:Sh} V.~Periwal and D.~Shevitz, {\it Unitary-Matrix Models as Exactly
Solvable String Theories}, Phys. Rev. Lett. {\bf 64} (1990), 1326--1329.
%
\bibitem{Sa:To} E. B. Saff and V. Totik, Logarithmic Potentials With External 
Fields, Springer-Verlag, New York, 1997. 
%
\bibitem{Sc} I. Schur, {\it Affektlose Gleichungen in der Theorie der 
Laguerreschen und Hermiteschen Polynomes}, J. f\"ur die reine und angewandte 
Matematik, {\bf 165} (1931), 52--58.
%
\bibitem{Se} A. Selberg, Bemerkninger om et multiplet integral, Norsk Mat. 
Tidsskr. {\bf 26} (1944), 71--78. 
%
\bibitem{St1} T. J.  Stieltjes, {\it Sur quelques th\'eor\`emes d'alg\`ebre}, 
Comptes Rendus de l'Academie des Sciences, Paris {\bf 100} (1885), 439--440. 
Oeuvres Compl\`etes, volume 1, 440--441.
%
\bibitem{St2} T. J. Stieltjes, {\it Sur les polyn\^omes de Jacobi}, 
Comptes Rendus de l'Academie des Sciences, Paris {\bf 100} (1885), 620--622. 
Oeuvres Compl\`etes, volume 1, 442--444.
%
\bibitem{Sz1} G.~Szeg\H{o}, {\it Ein Beitrag zur Theorie der Thetafunktionen}, 
Sitz. Preuss. Akad. Wiss. Phys. Math. Kl., {\bf XIX} (1926), 242--252,
Reprinted in "Collected Papers",  edited by R. Askey, Volume I, Birkhauser, 
Boston, 1982.
%
\bibitem{Sz2}  G. Szeg\H{o}, Orthogonal Polynomials, Fourth Edition, 
Amer. Math. Soc., Providence, 1975. 
%
\bibitem{Tr:Wi}  C. A. Tracy and H. Widom, {\it Random {U}nitary {M}atrices,
{P}ermutations and {P}ainl\'eve}, Commun. Math. Phys., {\bf 207}(1999),
665--685. 
%
\bibitem{Wi:Forr}  N. S. Witte and P. J. Forrester, 
{\it Gap probabilities in the finite and scaled {C}auchy random matrix 
ensembles}, in preparation, 2000 
%
\end{thebibliography}
\end{document}